\documentclass[]{interact}

\usepackage{epstopdf}
\usepackage[caption=false]{subfig}
\usepackage[doublespacing]{setspace}

\usepackage[natbibapa,nodoi]{apacite}
\usepackage{amssymb}
\usepackage{amstext}
\usepackage{amsmath}
\newcommand{\norm}[1]{\left\lVert#1\right\rVert}
\usepackage{amsthm}
\usepackage{algorithm}
\usepackage{algorithmic}
\usepackage{threeparttable}
\usepackage{mathtools}
\DeclarePairedDelimiter\ceil{\lceil}{\rceil}

\usepackage{array}
\newcolumntype{C}[1]{>{\centering\let\newline\\\arraybackslash\hspace{0pt}}m{#1}}
\newcolumntype{L}[1]{>{\left\let\newline\\\arraybackslash\hspace{0pt}}m{#1}}
\usepackage{multicol}
\usepackage{multirow}
\usepackage{bigstrut}

\theoremstyle{plain}
\newtheorem{theorem}{Theorem}[section]

\theoremstyle{definition}

\theoremstyle{remark}

\begin{document}


\title{A Note on the Stochastic Ruler Method for Discrete Simulation Optimization}

\author{
\name{Varun Ramamohan\textsuperscript{a}\thanks{CONTACT Varun Ramamohan. Email: varunr@mech.iitd.ac.in}, Utkarsh Agrawal\textsuperscript{a}, Mohit Goyal\textsuperscript{a}}
\affil{\textsuperscript{a}Department of Mechanical Engineering, Indian Institute of Technology Delhi, Hauz Khas, New Delhi 110016, India}
}

\maketitle

\begin{abstract}
In this paper, we propose a relaxation to the stochastic ruler method originally described by Yan and Mukai in 1992 for asymptotically determining the global optima of discrete simulation optimization problems. The `original' version of the method and its variants require that a candidate solution passes a certain number of tests with respect to the stochastic ruler to be selected as the next estimate of the optimal solution. This requirement can lead to promising candidate solutions being rejected and thus hinder the convergence of the method. Our proposed modification to the stochastic ruler method relaxes this requirement, and we show analytically that our relaxation incurs smaller computational overheads when a new solution in the neighborhood of the current solution is a `successful' candidate for the next estimate of the current solution. We also provide the theoretical grounding for the asymptotic convergence in probability of the variant to the global optimal solution under the same set of assumptions as those underlying the original stochastic ruler method. We then show numerically that this can yield improved convergence behavior - in terms of reducing the proportion of times the method does not achieve the specified performance criterion - via multiple numerical examples.
\end{abstract}

\begin{keywords}
Stochastic ruler; Discrete simulation optimization; Random search; Calibrator selection; Stochastic facility location
\end{keywords}

\section{Introduction}
\label{sec:intro}

Many problems in engineering and management science can be formulated as finding the best alternative, in terms of minimizing or maximizing an objective function, among a discrete, finite set of alternatives. Well-established and reasonably efficient methods for solving such optimization problems exist when the problems are deterministic; however, solving such problems with stochastic elements is a challenging task, especially when problem sizes are large, or when the set of alternatives are not easily enumerated \citep{henderson2006}. Such discrete stochastic optimization problems can be expressed as follows:

\begin{equation}
	\label{eq1}
	\begin{aligned}
		\underset{x \in \mathbb{S}}{\text{min}} \; f(x) = E[h(x, Y_x)]
	\end{aligned}
\end{equation} 

Here $\mathbb{S}$ is a discrete set, $h$ is a deterministic real-valued function and $Y_x$ is a random variable with its distribution dependent on the parameter $x \in \mathbb{S}$. Note that $h(x, Y_x)$ can itself be written as a random variable $H(x)$, in which case we can write $f(x) = E[H(x)] \; \forall \; x \in \mathbb{S}$.  Now let $S^* = \{x^* \in \mathbb{S} | f(x^*) \leq f(x) \; \forall \; x \in \mathbb{S}\}$. Thus $S^*$ represents the set of global optimal solutions to problem~\ref{eq1}, written as $S^* = \{x \in \mathbb{S} \; | \; f(x) \leq f(x') \; \forall \; x' \in \mathbb{S} \}$.

In many cases, closed form approximations for $f(x)$ or even $h(x, Y_x)$ may not be available. In such situations, the objective function in formulation~\ref{eq1} may have to be estimated via simulation (e.g., via Monte Carlo or discrete-event simulation). Thus the problem becomes a discrete simulation optimization problem. A variety of methods have been suggested to solve such problems, ranging from ranking and selection techniques \citep{jian2015intro} to random search methods \citep{andradottir2006}. The stochastic ruler method is one such discrete stochastic optimization method that has been proven to have asymptotic convergence in probability to the global optimum solution. The method was first developed by \citet{yan1992}, after which various modifications were published \citep{alrefaei1996dis,alrefaei2001mod,alrefaei2005}. In this paper, we present a simple relaxation of the original stochastic ruler method developed by \citet{yan1992} that may incur lesser computational expense under certain conditions, while still fitting within the theoretical framework required for convergence to the global optimum. We demonstrate the performance of the relaxation in comparison to the original method with four numerical examples. 

The stochastic ruler method as developed by \citet{yan1992}, and also with respect to its variants, involves starting with an initial solution and exploring its neighborhood (constructed to satisfy certain assumptions of reachability) to find a better solution. The method derives its name from the fact that a candidate solution from the previous solution's neighborhood is chosen by comparing each replicate value of the objective function associated with a candidate solution (generated from a stochastic simulation) to a sample from a uniform random variable $\theta(a,b)$ with parameters $a$ and $b$ ($a < b$) chosen to encompass possible values of $f(x)$. This uniform random variable represents the stochastic ruler. At the $k^{th}$ iteration of the method (where an iteration is defined as an evaluation of a single candidate solution), a candidate solution $x_{c(k)}$ is chosen as the next estimate $x_k$ of the optimal solution $x^*$ if it ``passes" every one of $M_k$ tests involving the stochastic ruler. Per the original stochastic ruler method, $M_k$ is chosen as an increasing function of $k$, and thus as the number of iterations ($k$) increase, passing $M_k$ tests also becomes less probable, and hence the algorithm is liable to becoming trapped within the neighbourhood of the solution estimate at an iteration. We propose a simple relaxation of this method in this paper: that a new candidate solution can be chosen as the next estimate of the optimal solution $x^*$ if it passes $\lceil{\alpha M_k}\rceil$ ($0 < \alpha < 1$) tests, $0 < \alpha < 1$. We show in the following sections: (a) that the relaxation fits within the probabilistic framework for asymptotic convergence of the original stochastic ruler method constructed by \citet{yan1992}, (b) that for `successful' candidate solutions, fewer tests need to be conducted for a given candidate solution than the original algorithm, and (c) via four numerical examples, that this simple relaxation can incur lesser computational expense and/or move towards the optimal solution with smaller probabilities of getting trapped in a given neighbourhood when compared to the original version of the stochastic ruler method. 

We now provide a brief overview of the literature. We do not provide a comprehensive overview of the many methods proposed to solve problems of the form in formulation~\ref{eq1} as it is beyond the scope of this paper; however, we discuss methods directly related to stochastic ruler method as used for simulation optimization. For a comprehensive accounting of discrete stochastic optimization methods, we direct the reader to relatively recent reviews \citep{amaran2016,de2019discrete}. 

As mentioned before, the stochastic ruler method was developed by \citet{yan1992}. In this seminal paper, they described in detail the stochastic framework for the asymptotic convergence in probability to the global optimum. Note that a key difference between this method and simulated annealing \citep{aarts1989simulated} is that the objective function at each iteration in simulated annealing is compared to its previous values whereas the uniform random variable is used in the stochastic ruler method. Note that simulated annealing has also been proposed as a method to solve the discrete stochastic optimization method \citep{gelfand1989}. 

Alrefaei and Andradottir \citep{alrefaei1996dis,alrefaei2001mod,alrefaei2005} develop modified versions of the stochastic ruler method wherein the algorithm generates a sequence of states (estimates of the solution), and (a) in one variant \citep{alrefaei2001mod}, the state (solution) that is visited most often is taken as the estimate of the optimal solution, and (b) in another variant \citep{alrefaei2005}, the state with the best average objective function value among all states visited is taken as the optimal solution. As part of these methods, a constant (and not increasing, as in the original stochastic ruler method) number of comparisons of the objective function against the stochastic ruler are performed. Even in these modified versions of the stochastic ruler method, all tests against the stochastic ruler must be successful for selection of the candidate solution as the new estimate of the optimal solution (i.e., a transition between states occurs). These modifications appear to build upon previous work by \citet{andradottir1995method} wherein the author proposes a discrete stochastic optimization method that utilizes the same optimality criterion. Note that in \citet{alrefaei2001mod} and \citet{alrefaei2005}, the authors also demonstrate how their approach can be used with a simulation in both transient and steady states, and also that their variants converge to the global optimum almost surely, as opposed to the in probability convergence of Yan and Mukai's method. We also refer the reader here to related work on random search \citep{andradottir1996global,andradottir1999,gong2000} and stochastic approximation methods for discrete stochastic optimization \citep{kleywegt2002}.

There appears to be limited work on the stochastic ruler method and its variants other than the literature discussed above. Most discrete simulation optimization research in recent years appears to have focused on achieving a balance between exploration, exploitation and estimation, and in locally convergent random search methods such as COMPASS. Instances of the former attempt to strike a balance between searching the solution space for global optima while at the same time extracting as much information locally as possible \citep{andradottir2009,hu2007}. The latter approach involves development of a framework for algorithms that provably converge to local optima \citep{hong2006,hong2007}.

In relation to the above literature, we provide a simple relaxation of the stochastic ruler method that decreases the probability of getting trapped in a given state (i.e., at a particular estimate of the optimal solution), especially when the solution being evaluated is in actuality a `successful' solution. In particular, this is likely to occur in the latter stages of the execution of the algorithm, as the number of tests that a candidate solution must pass increase, and as the solution being evaluated is not far from the global optimum. We focus on the stochastic ruler method because of its simplicity of implementation for many practical problems - in fact, we developed this relaxation while attempting to apply the method to the problem of selecting optimal calibration standards for a clinical measurement process \citep{wsc18vr}. Thus, as part of this problem, we also describe the application of the relaxation to multiple detailed numerical examples. This may be useful to practitioners given that a key stage in applying this method to a problem involves defining the neighborhood structure for candidate solutions, and we anticipate that our examples will add to the relatively small set of examples in the literature for the implementation of this method. 

In Section~\ref{algo}, we describe the original stochastic ruler method and our proposed relaxation, and how our proposed relaxation fits within the global convergence framework developed by \citet{yan1992}. In Section 3, we describe the application of our proposed relaxation to four numerical examples, and benchmark its performance against the original stochastic ruler method. We conclude in Section 4 with a summary of our work, its limitations and avenues for ongoing and future research.

\section{Stochastic Ruler: Original Algorithm and Proposed Relaxation}
\label{algo}
We begin by providing a brief description of the original stochastic ruler algorithm, and then describe our proposed relaxation in relation to the original.

\subsection{Stochastic Ruler Algorithm}
\label{orsr}
The stochastic ruler (SR) method involves conducting a number of tests to evaluate a given solution $x \in \mathbb{S}$, wherein each test involves comparing a sample of the $h(x, Y_x)$ with a sample from $\theta(a,b)$. Given that formulation~\ref{eq1} is a minimization problem, the test is considered successful if $h(x, Y_x) \leq \theta(a, b)$. Note that each $h(x, Y_x)$ is generated from a simulation (e.g., a Monte Carlo simulation). Thus the minimization problem in formulation~\ref{eq1} is converted to a stochastic maximization problem of the form:

\begin{equation}
	\label{eq2}
	\max \; Pr(h(x, Y_x) \leq \theta(a,b) | \forall \; x \in \mathbb{S})
\end{equation}

The SR method attempts to solve this maximization problem in place of formulation~\ref{eq1}. Note that we henceforth write the probability $Pr(h(x, Y_x) \leq \theta(a,b))$ in short as $P(x, a, b)$, implying that formulation~\ref{eq2} can be written more briefly as $\max\;\{P(x, a, b) | x \in \mathbb{S}\}$. The optimal set of solutions to the problem in formulation~\ref{eq2} thus can be written as:

\begin{equation}
	\nonumber
	S^*(a,b) = \{x \in \mathbb{S} \; | \; P(x,a,b) \geq P(x',a,b) \; \forall \; x' \in \mathbb{S} \}
\end{equation}

We now provide in brief the key definitions and assumptions required by the method.

\textbf{Definition 1} from \cite{alrefaei1996dis}. For each $x \in \mathbb{S}$, there exists a subset $N(x) \subseteq \mathbb{S} - x $ which is called the set of neighbors of $x$.

\textbf{Definition 2} from \cite{alrefaei1996dis}. A function $R : \mathbb{S} \times \mathbb{S} \to [0, 1]$ is said to be a transition probability for $\mathbb{S}$ and $N$ if:	$R(x, x') > 0 \implies x' \in N(x)$, and $\sum_{x' \in \mathbb{S}} R(x, x') = 1, \forall \; x \in \mathbb{S}$. 

Further, the SR method demands that the $H(x) = h(x, Y_x)$ has finite variance; that is, $E(H(x)^2) < \infty~, \forall~ x \; \in \mathbb{S}$. 

We also note here that $P(x,a,b)$ essentially represents the probability that a single test conducted as part of evaluating a candidate solution $x \in \mathbb{S}$ is successful; thus, the probability that $x$ is selected as the new estimate of the solution is a function of $P(x,a,b)$. Therefore, given that our relaxation changes the probability that a candidate solution is selected as the new estimate of the solution, any condition that the original method imposes on $P(x,a,b)$ must also be satisfied by the version of $P(x,a,b)$ used in our relaxation of the method.

The relationship between the maximization problem in formulation~\ref{eq2} and the minimization problem in formulation~\ref{eq1} is given by Yan and Mukai in the following theorem. 

\begin{theorem}
	\label{thm1}
	Theorem 3.1 from \cite{yan1992}. There exist real numbers $a'$ and $b'$ such that $a' < b'$ and for any $a < a'$ and any $b > b'$, the following conclusions can be made:
	
	(1) If $f(x) < f(x')$ then $P(x, a, b) > P(x', a, b)$, 
	
	(2) $0 < P(x, a, b) < 1$, $\forall \; x \in \mathbb{S}$,
	
	(3) $S^*(a, b)$ $\subset$ $S^*$ and $S^*(a, b) \neq \phi$.
	
\end{theorem}

Thus if our relaxation is to fit within the theoretical framework developed by Yan and Mukai, then our versions of $P(x,a,b)$ must also satisfy the conditions in the theorem above.

The following assumptions must also be made for the method. These correspond to Assumptions 2 - 5 from \cite{alrefaei1996dis}.

\textbf{Assumption 1}. For any $x, x' \in \mathbb{S}$, $x'$ is reachable from $x$. That is, there exists a finite sequence $\{n_i\}_{i = 0}^j$ for some $j$, such that $x_{n_0} = x, \; x_{n_j} = x'$, and $x_{n_{i+1}} \in N(x_{n_i}), \; i = 0, 1, 2, . . , j-1$.

\textbf{Assumption 2}. The neighbour system $N$ and the transition probability $R$ are symmetric. That is,

1. $x' \in N(x) \implies x \in N(x')$, and

2. $R(x, x') = R(x', x) \; \forall \; x,\; x' \in \mathbb{S}$.

\textbf{Assumption 3}. The parameters $a, b \in \mathbb{R}$ are selected to cover the range of the observed objective function values $H(x)$, where $x \in \mathbb{S}$.

\textbf{Assumption 4}. A sequence $\{M_k\}$ of positive integers satisfies $M_k \to \infty$ as $k \to \infty$.

We are now in a position to provide Yan and Mukai's SR algorithm. 

\noindent \textbf{Original SR algorithm}. 

\textit{Initialization}: neighborhood structure $N$ and transition probability function $R$ for $x \in \mathbb{S}$, monotonically increasing sequence $\{M_k\}$, SR $\theta(a,b)$, initial solution $x_0$, and set $k = 0$.

\textit{Step 1}. Given $x_k = x$, choose a new candidate solution $z_k$ from the neighborhood $N(x)$ with probability distribution $P\{z_k = z \; | \; x_k = x\} = R(x, z), \; z \in N(x)$.

\textit{Step 2}. Given $z_k = z$, generate a realization $h(z)$ of $H(z)$ (i.e., generate one replicate value from the simulation). Then generate a realization $\theta$ from $\theta(a, b)$. If $h(z) > \theta$, then let $x_{k+1} = x_k$ and go to Step 3. If not, generate another realization $h(z)$ from $H(z)$ and another realization $\theta$ from $\theta(a, b)$. If $h(z) > \theta$, then let $x_{k+1} = x_k$ and go to Step 3. Otherwise continue to generate realizations and conduct the tests. If all $M_k$ tests, $h(z) > \theta$, fail, then accept the candidate $z_k$ and set $x_{k+1} = z_k = z$.

\textit{Step 3}. Set $k = k + 1$ and go to Step 1.  

Conducting the $M_k$ tests in \textit{Step 2} above contributes the majority of the computational expense of implementing the method, in particular for successful candidate solutions. Further, we note that each time one test fails, the number of tests to be conducted for another candidate solution in the same neighborhood increases by one. This implies that the number of tests that need to be conducted to evaluate a solution can quickly increase, depending upon the choice of sequence $\{M_k\}$.

Note that the stochastic process $\{x_k\}$ produced by the algorithm above is a discrete-time Markov chain defined over states $x \in \mathbb{S}$, and the associated state transition probabilities are given by:

\begin{equation}
	\label{eq3}
	\begin{aligned}
		&P_{x x'} (M_k) = P[x_{k+1} = x' \; | \; x_k = x] \; ; \\
		& = R(x,x') P(x',a,b)^{M_k} \; \text{if} \;\; x' \in N(x) \; ; \\
		& = 1 - \sum_{x' \in \mathbb{S}} R(x,x') P(x',a,b)^{M_k} \;\;\text{if} \; x' = x \; ; \\
		& = 0 \; \text{otherwise}
	\end{aligned}
\end{equation} 

The algorithm is typically terminated once a sufficient decrease in the objective function value is achieved, or the computational budget is spent. We also note here that if the number of tests $M_k$ is set to a constant number $M$ in each iteration, the resulting Markov chain becomes stationary. 

The convergence framework for the asymptotic (in probability) convergence of the SR algorithm to the global optimum consists of the following steps. 
\begin{itemize}
	\item[1.] After the correspondence between formulation~\ref{eq1} and formulation~\ref{eq2} is established via Theorem~\ref{thm1}, the authors set the parameter $M_k$ equal to a constant $M$, and analyze the convergence of the resulting stationary Markov chain to its stationary distribution $\pi(M)$.
	\item[2.] They then analyze the behaviour of $\pi(M)$ as $M \to \infty$.
	\item[3.] Next, they consider the behaviour of the nonstationary Markov chain when $M_k$ is not a constant, and prove that this Markov chain is strongly ergodic. This leads to the result that the $\{x_k\}$ lie in the set of optimal solutions $S^*$ as $k \to \infty$.  
\end{itemize} 

Note that Yan and Mukai also consider the rate of convergence of the SR algorithm in a subsequent step; however, this is not necessary for proving convergence to a global optimum, and hence we do not discuss this in this paper.

Our approach involves a change in the transition probabilities $P_{x x'}$ associated with the Markov chain. Therefore, we describe steps 1 and 2 in the framework above in more detail so that we can demonstrate how the results in these steps also apply to our approach. It will then follow that the results in step 3 also will apply to our approach. 

In step 1 above, Yan and Mukai first consider $M_k$ to be a constant $M$, and define a vector $\pi(M)$ of dimension $|\mathbb{S}|$ with elements given by:

\begin{equation*}
	\pi_x(M) = \frac{\{P(x,a,b)\}^{M}}{\sum_{x' \in \mathbb{S}}\{P(x',a,b)\}^M} 
\end{equation*}

In the below theorem, the authors then prove that the vector defined above is the stationary distribution for the Markov chain defined by the elements of $\mathbb{S}$.

\begin{theorem}
	\label{thm2}
	\textit{Theorem 5.1 from \cite{yan1992}.} The vector $\pi(M)$ with elements $\pi_x(M)$ as defined above is the stationary probability distribution for the Markov chain $\{x_k\}$ defined by transition probabilities in equation~\ref{eq3}, given that $M_k = M$.
\end{theorem}

In step 2 of the convergence framework, Yan and Mukai begin their analysis of the behaviour of $\pi(M)$ as $M\to \infty$ by defining a probability vector.

\textbf{Definition 3} (Definition 6.1 from \cite{yan1992}). For a finite set $\mathbb{S}$, the set $\prod(\mathbb{S})$ of positive unit vectors is called the set of probability vectors for $\mathbb{S}$, if:
\begin{equation*}
	\Pi(S) = \{\pi \in [0,1]^{\kappa} \; | \; \pi_x \geq 0, \norm{\pi} = \sum_{x \in \mathbb{S}} \pi_x  = 1\},
\end{equation*}

where $\kappa = |\mathbb{S}|$.

The authors then define the notion of an optimal probability vector as follows: a probability vector $\pi^*$ for $\mathbb{S}$ is referred to as \textit{optimal} if $\pi^*_x = 0$ for any $x \notin S^*$. Further, if $\pi_x^* = 0$ for any $x \notin S^*(a,b)$, then $\pi^*$ is optimal and $x \in S^*$ given that $S^*(a,b) \subset S^*$.

The following result is then proved regarding the asymptotic behaviour of $\pi(M)$ as $M \to \infty$. 
\begin{theorem}
	\label{thm3}
	\textit{Theorem 6.1 from \cite{yan1992}.} The probability vector $\pi(M)$ defined above converges, as $M \to \infty$, to an optimal probability vector $\pi^*$. Furthermore, we have: 
	\begin{equation*} 
		\begin{aligned}
			&\pi^*_x = 1 / |S^*(a, b)|, \; \text{if} \; x \in S^*(a, b), \\
			&\pi^*_x = 0, \; \; \text{otherwise},
		\end{aligned}
	\end{equation*}
	where $|S^*(a, b)|$ represents the cardinality of $S^*(a, b)$.
\end{theorem}

Thus, given that our relaxation involves changes to $P_{xx'}$ and thus to $\pi(M)$, we need to prove that versions of the above results can be proved for the $P_{xx'}$ and $\pi(M)$ that we generate for our approach.

We are now in a position to describe in detail our proposed relaxation of the SR method.

\subsection{Stochastic Ruler Algorithm: Proposed Relaxation}
\label{srrelax}
In our proposed relaxation to the original SR algorithm, we require only a fraction $\alpha,\; 0 < \alpha < 1$, of $M_k$ tests to be successful for the candidate solution to be chosen as the new estimate of the optimal solution. This implies that in comparison to the original SR method, lesser computational effort is likely to be expended to evaluate candidate solutions that are genuinely better than the current estimate, and conversely, that there is a larger probability (depending upon how close $\alpha$ is to 1) of an inferior candidate being chosen as the new estimate of the optimal estimate.

Note that we make use of the definitions and assumptions regarding the neighborhood structure associated with the original method for this algorithm as well. The method requires the following steps.\\ 

\noindent \textbf{Relaxation of SR algorithm}. 

\textit{Initialization}: neighborhood structure $N$ and transition probability function $R$ for $x \in \mathbb{S}$, monotonically increasing sequence $\{M_k\}$, stochastic ruler $\theta(a,b)$, initial solution $x_0$, $\alpha\; (0 < \alpha < 1)$, and set $k = 0$.

\textit{Step 1}. Given $x_k = x$, choose a new candidate solution $z_k$ from the neighborhood $N(x)$ with probability $P\{z_k = z \; | \; x_k = x\} = R(x, z), \; z \in N(x)$.

\textit{Step 2}. Given $z_k = z$, generate a realization $h(z)$ of $H(z)$ (i.e., generate one replicate value from the simulation). Then generate a realization $\theta$ from $\theta(a, b)$. Check whether $h(z) > \theta$. If not, increment the number of ``successful" tests by one. If the number of successful tests equals $\lceil\alpha M_k\rceil$, accept the candidate $z_k$ and set $x_{k+1} = z_k = z$. If the number of unsuccessful tests exceeds $M_k - \lceil\alpha M_k\rceil$, then let $x_{k+1} = x_k$ and go to Step 3.


\textit{Step 3}. Set $k = k + 1$ and go to Step 1. 

It is evident that the relaxation is only in Step 2 above. Thus the state transition probabilities for the Markov chain $\{x_k\}$ generated by our relaxation of the SR method also change as given below. For the below expression, we let $n = \lceil\alpha M_k\rceil$.

\begin{equation}
	\label{eq4}
	\begin{aligned}
		&P_{x x'} (M_k) = P[x_{k+1} = x' \; | \; x_k = x] \; ; \\
		& = R(x,x') \left[\sum_{t = n}^{M_k} \; {t-1 \choose t-n} \; P(x',a,b)^n  (1- P(x',a,b))^{t-n} \right] \; \text{if} \;\; x' \in N(x) \; ; \\
		& = 1 - \sum_{x' \in \mathbb{S}} R(x,x') \left[\sum_{t = n}^{M_k} \; {t-1 \choose t-n} \; P(x',a,b)^n  (1- P(x',a,b))^{t-n} \right] \;\;\text{if} \; x' = x \; ; \\
		& = 0 \; \text{otherwise.}
	\end{aligned}
\end{equation}

The above expression for $P_{xx'}(M_k), \; \text{when} \; x' \in N(x)$, follows from the fact that if $t \; (t \geq n)$ tests are conducted and $x'$ is selected as the new estimate of the optimal solution from these $t$ tests, the last test conducted must be a success. Further, we also see that we use the same definition of $P(x, a, b)$ as the original SR method, so Theorem~\ref{thm1} applies to our method as well. 

In order to demonstrate that the convergence of our relaxation of the SR method can be accommodated within the framework of the original method, we also define a probability vector considering $M_k$ constant and equal to some $M$. We refer to this probability vector as $\pi_r(M)$, with each of its elements given by the expression below (letting $\lceil\alpha M\rceil = n$).

\begin{equation}
	\label{eq5}
	\pi_{x(r)}(M) = \frac{\sum_{t = n}^{M} \; {t-1 \choose t-n} \; P(x,a,b)^n  (1- P(x,a,b))^{t-n}}{\sum_{x' \in \mathbb{S}}\left[\sum_{t = n}^{M} \; {t-1 \choose t-n} \; P(x',a,b)^n  (1- P(x',a,b))^{t-n} \right]} 
\end{equation}

We now state the equivalent of Theorem~\ref{thm2} for $\pi_{r}(M)$. 

\begin{theorem}
	\label{thm4}
	The vector $\pi_{r}(M)$ with elements $\pi_{x(r)}(M)$ as defined in equation~\ref{eq5} is the stationary probability distribution for the Markov chain $\{x_k\}$ defined by transition probabilities in equation~\ref{eq4}, given that $M_k = M$.
\end{theorem}

\begin{proof}
	The proof of the above result is similar to that given in the Appendix of \cite{yan1992} for Theorem~\ref{thm2}, with the terms from equations~\ref{eq4} and~\ref{eq5} replacing the corresponding terms in the proof by Yan and Mukai.
\end{proof}

The next step is to state and prove the equivalent of Theorem~\ref{thm3} for the proposed relaxation. Note that we use the same definition of an optimal probability vector $\pi_r^*(M)$ for our relaxation of the SR method as well. 

\begin{theorem}
	\label{thm5}
	The probability vector $\pi_r(M)$ converges, as $M \to \infty$, to an optimal probability vector $\pi_r^*$. Furthermore, we have: 
	\begin{equation*} 
		\begin{aligned}
			&\pi^*_{x(r)} = 1 / |S^*(a, b)|, \; \text{if} \; x \in S^*(a, b), \\
			&\pi^*_{x(r)} = 0, \; \; \text{otherwise},
		\end{aligned}
	\end{equation*}
	where $|S^*(a, b)|$ represents the cardinality of $S^*(a, b)$.
\end{theorem}

The proof for the above theorem, as provided in Yan and Mukai, follows from the definition of the stationary distribution $\pi_r(M)$ in equation~\ref{eq5} and from Theorem~\ref{thm1}.

The rest of the proof of the asymptotic convergence of the proposed relaxation of the SR method to the global optima follows the same steps as the paper in \cite{yan1992} (i.e., step 3 of the convergence framework).

\subsection{Expected Number of Tests for Successful Candidate Solutions}

Before we provide numerical results, we consider analytically the performance of the proposed relaxation of the SR method in comparison with the original version with respect to the key improvement that we suggest results from the relaxation: that less computational effort will be expended with the relaxation when a candidate solution in the neighborhood of the current estimate of the solution is indeed an acceptable next estimate of the solution. In order to establish this analytically, we would have to prove that the expected number of comparisons between samples from $H(z)$ and $\theta(a, b)$ would be lesser than the number of comparisons ($M_k$, at the $k^{th}$ candidate solution evaluated) required with the original version of the SR method, given that the candidate solution evaluated is indeed selected as the new estimate of the current solution. While this would appear to be an intuitively obvious result, proving this under the framework of the relaxation of the SR method helps establish analytically the principal advantage of the relaxation in comparison to the original version. 

We begin by defining the term $s$ as the complement of $P(x, a, b)$. We remind readers here that $P(x,a,b)$ is the probability that a given `test' or comparison is successful; that is, $P(x,a,b) = P(h(x) \leq \theta)$, where $h(x)$ is a single replicate value of the simulation (i.e., a sample from $H(x)$), and $\theta$ is a sample from the stochastic ruler $\theta(a,b)$. In other words, $s = 1 - P(x,a,b)$. If $M$ is the maximum number of comparisons that can be performed when evaluating a particular solution $x$, then we can write for the original SR method:

\begin{equation*}
	\begin{aligned}
		&P(\text{$x$ is selected as new solution estimate}) = (1 - s)^M\\
		&P(\text{$x$ is not selected as new solution estimate}) = 1 - (1 - s)^M
	\end{aligned}
\end{equation*}

If $x$ is selected as the new estimate of the optimal solution, then under the original SR method, the number of comparisons performed is always $M$. Under the relaxation of the SR method, we define a random variable $T$ (with realizations $t$) that represents the number of tests performed when evaluating a given candidate solution $x$. If $x$ is chosen as the new estimate of the solution, then $T$ will take values between $n = \ceil{\alpha M}$ and $M$. We now define two events $A$ and $B$ as follows: $A$ is the event wherein $x$ is chosen as the new estimate of the optimal solution when exactly $t$ tests ($t > n$) are conducted, and $B$ is the event wherein $x$ is chosen as the new estimate of the optimal solution when $t$ takes any value between $n$ and $M$ (inclusive). We are thus interested in proving that $E[A|B] < M$. 

We begin by noting that the probability of selecting $x$ as the new estimate of the optimal solution in exactly $t$ tests (i.e., the probability that the event $A$ occurs) can be written as follows.

\begin{equation*}
	\begin{aligned}
		&P(A) = {t-1 \choose t - n} s^{(t-n)} (1 - s)^n, \; t \geq n\\
		&     \;\;\;\;\;\;\;\;\;=  0, \;t < n
	\end{aligned}
\end{equation*}

The above expression follows from the fact that if event $A$ occurs upon performing $t$ tests, then the last ($t^{th}$) test must be successful. We also note that the probability of event $B$ can similarly be written as follows.

\begin{equation*}
	P(B) = \sum\limits_{t = n}^{M}\; {t-1 \choose t - n} s^{(t-n)} (1 - s)^n, \; t \geq n
\end{equation*}

We are interested in determining $P(A|B)$ for calculating $E[A|B]$, and for this we utilize Bayes' Theorem as follows. 

\begin{equation*}
	P(A|B) = P(B|A) \times P(A) / P(B)
\end{equation*}

Here, we note that $P(B|A) = 1$, by the definition of the events $A$ and $B$. Therefore, we can write:

\begin{equation}
	\label{eqagivb}
	P(A|B) = \frac{{t-1 \choose t - n} s^{(t-n)} (1 - s)^n}{\sum\limits_{t = n}^{M}\; {t-1 \choose t - n} s^{(t-n)} (1 - s)^n}
\end{equation}

Note that equation~\ref{eqagivb} is nonzero for only $t \geq n$ based on the definitions of $A$ and $B$. Now $E[A|B]$ can be written as follows.

\begin{equation}
	\begin{aligned}
		&E[A|B] = \sum\limits_{t = n}^M t \times P(A|B)\\
		&E[A|B] = \frac{1}{P(B)}\;\sum\limits_{t = n}^M t \times {t-1 \choose t - n} s^{(t-n)} (1 - s)^n
	\end{aligned}
\end{equation}

Now, we can write:

\begin{equation}
	\label{eqm}
	\begin{aligned}
		&E[A|B] = \frac{1}{P(B)}\;\sum\limits_{t = n}^M t \times {t-1 \choose t - n} s^{(t-n)} (1 - s)^n < \frac{1}{P(B)}\;\sum\limits_{t = n}^M M \times {t-1 \choose t - n} s^{(t-n)} (1 - s)^n
	\end{aligned}
\end{equation}

Now we note that $P(B) = \sum\limits_{t = n}^{M}\; {t-1 \choose t - n} s^{(t-n)} (1 - s)^n$, and hence we can see from equation~\ref{eqm} that:

\begin{equation}
	\label{eqm}
	E[A|B] < \frac{M}{P(B)} \times P(B) \; \implies \; E[A|B] < M
\end{equation}

Thus we see that if a candidate solution is to be chosen as the next estimate of the optimal solution (a `successful' candidate solution), our proposed relaxation does so in fewer steps. On the other hand, if a candidate solution is to be rejected as the next estimate of the optimal solution (an `unsuccessful' candidate solution), then the original version of the SR method will do so in fewer steps. Thus, whether the average computational expense of our proposed relaxation of the SR method will be lesser than or greater than that of the original version depends on the distribution of successful and unsuccessful candidate solutions in the solution space, and the variances of these solutions. However, in situations where the proportion of successful candidates across the solution space is larger than that of unsuccesful candidates, and their variances are comparable, our relaxation of the SR method is likely to incur lesser computational expense.

We now provide numerical examples of the implementation of the proposed relaxation of the SR method and benchmark it with respect to the original SR method.

\section{Numerical Experiments}
\label{numsec}
\subsection{Numerical Examples with Known Global Minima}
We begin with the same numerical example provided in \cite{alrefaei1996dis}. 

\begin{figure}[!h]
	\vspace{-0.2cm}
	\centering
	\includegraphics[width=0.75\textwidth]{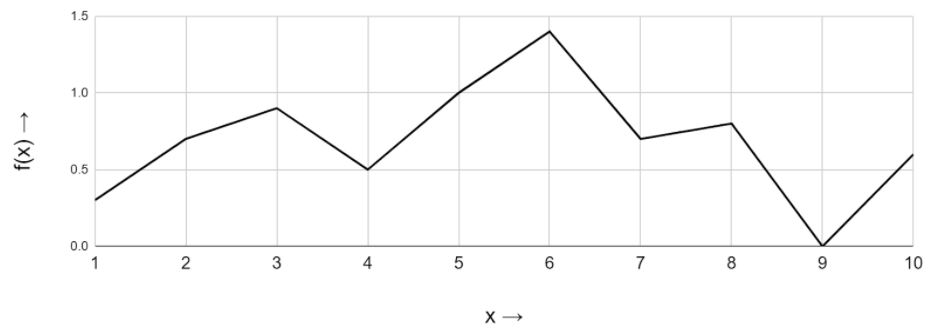}
	\caption{Objective function for the numerical example in \cite{alrefaei1996dis}.}
	\label{fig1}
	\vspace{-0.1cm}
\end{figure}

Let the search space $\mathbb{S} = \{1, ... ,10\}$, with $f(x) = E[Y_x] \; \forall \; x \in \mathbb{S}$. Here $Y_x$ is uniformly distributed on the interval $f(x) \pm 0.5,\; \forall \; x \in \mathbb{S}$, and $f(1), ... , f(10)$ are $0.3, 0.7, 0.9, 0.5$, $1.0, 1.4, 0.7, 0.8, 0.0, \text{and}\; 0.6$, respectively. We set $a = -0.5, b = 1.9$, and $x_0$ is sampled from $\mathbb{S}$ assuming a uniform distribution. Thus, the stochastic ruler for this example $\theta(a, b)$ is given by $\theta(-0.5, 1.9)$. The objective function for this problem is given in Figure 1 below (similar to Figure~\ref{fig1} from \cite{alrefaei1996dis}).

The neighbourhood structure used for this problem is given below. 

\begin{equation*}
	N(x) = S - \{x\}, \forall \; x \in \mathbb{S}.
\end{equation*}

Note that the global minimum for the objective function occurs at $x = 9$, whereas local minima also occur at $x = 7$ and $x = 4$. In Table~\ref{tab1}, we provide a comparison of the two algorithms on the basis of two metrics: (a) average runtime, and (ii) average number of iterations. An iteration here represents the evaluation of a single candidate solution; that is, a single execution of Step 2 of the original SR algorithm or our proposed relaxation. We terminate the algorithm when the global minimum is reached, given that we know the global minimum for this problem. The results in the table are the average of 500 replications of the implementation of the method.

\begin{table}[htbp]
	\centering
	\caption{Comparison between original and relaxation of the SR method: numerical example 1.}
	\begin{tabular}{|C{5cm}|C{3cm}|C{3cm}|}
		\hline
		Performance Metric & SR Original & SR Relaxation \\
		\hline
		Average runtime (seconds) & 2.26 $\times$ $10^{-4}$ & 2.468 $\times$ $10^{-4}$ \\
		\hline
		Average iterations & 11.244 & 11.244 \\
		\hline
	\end{tabular}%
	\label{tab1}%
\end{table}%

From the above results, while we see that the average runtime is marginally (approximately 9\%) lower for the original version, these results do not conclusively provide information regarding the performance of the relaxation of the SR method in comparison with the original version. This is likely due to the very small size of the problem.

Therefore, we constructed a larger version of the above problem, with 100 states instead of 10. For this problem, we let $\mathbb{S} = {1,2,3,...,100}$ and $f(x)$ correspondingly takes values in $D = \{0.3, 0.7, 0.9,..,1.9, 1.4\}$. Stochasticity in the problem is introduced by modeling $Y_x$ as a uniform random variable on $f(x) \pm 0.5$ (inclusive). Thus if we write the random variable that each replication of the simulation produces as $h(x, Y_x)$, the optimization problem becomes equivalent to problem~\ref{eq1}. The graph of the objective function is provided in Figure~\ref{numex2} below.

\begin{figure}[htbp]
	\vspace{-0.2cm}
	\centering
	\includegraphics[width=0.75\textwidth]{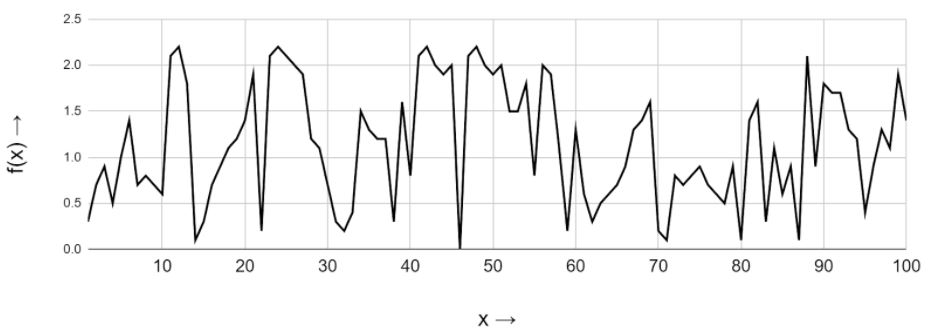}
	\caption{Objective function for numerical example 2.}
	\label{numex2}
	\vspace{-0.1cm}
\end{figure}

The problem has a single global minimum at $x = 46$. The limits of the stochastic ruler as taken as the lower and upper bounds of the uniform random variables associated with the minimum and maximum values of $x \in \mathbb{S}$; that is, we set $a = -0.5,\; b = 2.7$.

We develop a neighborhood stucture for this problem as follows. We assume each point $x$ has 10 neighbors, 5 of which are less than $x$ and 5 of which are greater than $x$. More precisely, the neighborhood $N(x)$ of $x$ is defined as $x \pm 1$, $x \pm 2$, $x \pm 3$, $x \pm 4$, $x \pm 5$. It is clear that if $x' \in N(x) > 100 \;\text{or}\; x' < 1$, then we leave the feasible region $\mathbb{S}$. In this situation, we simply subtract or add 100 to $x'$, respectively. We provide illustrative examples below. 

\begin{equation*}
	\begin{aligned}
		&N(43) = \{38, 39, 40, 41, 42, 44, 45, 46, 47, 48\}\\
		&N(3) = \{98, 99, 100, 1, 2, 4, 5, 6, 7, 8\} \\
		&N(98) = \{93, 94, 95, 96, 97, 99, 100, 1, 2, 3\}
	\end{aligned}
\end{equation*}

Further, $R(x,x'),  = 0.1 \; \forall \; x' \in N(x)$. The results of applying the original and the relaxation of the SR method are summarized in Table~\ref{tab2} below. The termination criterion is the same as that for the previous example: that is, the algorithm is terminated when the global minimum is found. 

\begin{table}[htbp]
	\centering
	\caption{Comparison between original and relaxation of the SR method: numerical example 2.}
	\begin{tabular}{|C{1cm}|C{2cm}|C{2cm}|C{2cm}|C{2cm}|}
		\hline
		\multirow{2}[4]{*}{$\alpha$} & \multicolumn{2}{c|}{SR Relaxation} & \multicolumn{2}{c|}{SR Original} \bigstrut\\
		\cline{2-5}   & Runtime (seconds) & Average $k$ & Runtime (seconds) & Average $k$ \bigstrut\\
		\hline
		0.75  & 0.0436 & 2568.088 & 0.1756 & 11704.88 \bigstrut[t]\\
		0.70   & 0.0254 & 1441.094 & 0.1756 & 11704.88 \\
		0.65  & 0.0172 & 898.806 & 0.1756 & 11704.88 \\
		0.60   & 0.0175 & 872.066 & 0.1756 & 11704.88 \\
		\hline
	\end{tabular}%
	\label{tab2}%
\end{table}%

The runtimes and the number of iterations (indicated by $k$) for the relaxation are much lower than those of the original SR method. Further, the performance seems to improve significantly as we lower the value of $\alpha$, likely because the feasible space is covered faster with decreasing $\alpha$. However, the value of $\alpha$ must be chosen with care, given that the probability of incorrect selection of an `unsuccessful' candidate solution as the next estimate of the optimal solution will increase with decreasing $\alpha$.

\subsection{Numerical Example 3: Optimal Assay Calibrator Selection}
\label{calib}

In both of the above numerical examples, we were aware of the global minimum prior to implementing the algorithm, and hence the termination criterion was also set accordingly. We now consider an example where the global minima are not known, and hence the termination criteria are specified on the basis of the available computational budget and whether a sufficient decrease in the objective function value is achieved. This numerical example is motivated by the work in \cite{wsc18vr}, which involved the selection of optimal calibrator concentrations for the calibration of a clinical measurement process such that the uncertainty of measurement around a test sample is minimized. Such clinical measurement processes include routine medical laboratory investigations such as determining a patient's cholesterol or blood glucose levels.

The measurement process we consider involves an indirect measurand - that is, the instrument measures some property of the sample, which is then converted into the quantity of interest via a calibration function. This relationship can be expressed as follows: $m = g(c)$, where $m$ is the quantity of interest (the measurand) and $c$ is the quantity actually measured by the instrument, and $g$, typically a map from $\mathbb{R}\to \mathbb{R}$, is the calibration function. The calibration process involves estimation of the parameters of the function $g$ using samples (say $n$ in number) with known values of the quantity of interest, which are referred to as calibrators, each with a known value of $m$ - for example a set of calibrators $m$, given by $m = \{m_1, m_2, ...,m_n\}$. Measurements of the set of calibrators $m$ are conducted to generate corresponding values of $c = \{c_1, c_2,...,c_n\}$ for each calibrator. These measurements $c$ and calibrators $m$ are then used to estimate the parameters of the calibration function. The calibrated function is then used to evaluate samples with unknown values of $m$ - that is, the actual measurement process. 

The measurement process can involve, for example, the measurement of $N$ samples - all with a desired concentration of $s$ (measurements conducted, for example, as part of a quality control process used to evaluate the measurement uncertainty of the measurement process) - and the conversion of the readings for these $N$ samples into the desired quantity by a calibration function parameterized as described above. This would yield $N$ different readings given by $S = \{s_1, s_2,...,s_N\}$ due to the stochasticity introduced by the sources of uncertainty operating within the measurement process. The measurement uncertainty, typically defined as the standard deviation of a set of measurements of samples with the same target measurand value, would then be estimated as $u_s(m) = \sqrt{\frac{\sum\limits_{i = 1}^N (s_i - \bar{s})^2}{N-1}}$. Here $u_s$ is written as an implicit function of $m$ given that the measurement uncertainty associated with a target concentration of $s$ would be dependent on the calibration function generated by choosing a particular set of calibrators $m$. Note that $u_s$ would be dependent on other factors as well, such as the distributions of the sources of uncertainty; however, the analyst may only have control over the choice of $m$, and hence we write $u_s$ as a function of only $m$.

The stochasticity in the calibration and unknown sample measurement processes comes from sources of uncertainty operating within these processes. Thus $c$ and $S$ are both random vectors. A detailed accounting of sources of uncertainty, their characterization, and incorporation into a mathematical model of measurement uncertainty for the haptoglobin clinical measurement process is given in \cite{wsc18vr}. This mathematical model of the measurement uncertainty is often not analytically tractable such that $u_s$ can be estimated manually, and thus Monte Carlo simulation is often employed in the estimation of $u_s$, as in the case of the haptoglobin assay. 

The haptoglobin assay measurement uncertainty model in \cite{wsc18vr} is based on a four point log-logit calibration function. However, because the purpose of this work is to demonstrate the application of our proposed relaxation of the SR method and not to describe the application itself in full detail, we do not work with an actual assay measurement uncertainty model. Instead, we consider a hypothetical assay with a quadratic calibration function, implying that three calibration function parameters are estimated as part of the calibration process. However, the construction of the model of measurement uncertainty of the assay, including the calibration process, is based on that for the haptoglobin assay in \cite{wsc18vr}. The parameters of the calibration function are estimated using six calibrator concentration values and their corresponding measurements via nonlinear regression, as is common practice for many assays with nonlinear calibration functions \citep{wsc15vr,wsc18vr}.

Uncertainty is incorporated into the calibration process via the calibrators and their corresponding measurements, which yields a calibration function with parameters different from those that would have been obtained had there been no measurement uncertainty. This perturbed calibration function is used to evaluate samples with unknown values of the measurand. The measurements of these samples are also subject to uncertainty. All the sources of uncertainty in the measurement process are assumed to be normally distributed with mean equal to the target value of each component (calibrator concentration or corresponding calibrator/sample measurements), implying no measurement bias, and a standard deviation equal to 5\% of the target value. The Monte Carlo simulation used to estimate the measurement uncertainty is depicted in Algorithm~\ref{alg1} below. Note that, as in \cite{wsc18vr}, a reference calibration function $l$ and its inverse are used to generate desired `error-free' values of the measurements (prior to the simulated introduction of measurement error from the sources of uncertainty within the measurement process) corresponding to the calibrator/sample value.

\begin{algorithm}
	\caption{Simulation model for estimating assay measurement uncertainty.}
	\label{alg1}
	\begin{algorithmic}
		\STATE \textit{Initialization:} Reference function $l$, number of replicate sample measurements $N$, test sample with target concentation $s_t$ and corresponding measurement $c_t  =l(s_t)$.\\
		/*begin calibration process*/\\
		\FOR{$j = 1-6$} 
		\STATE $m_j' \sim N(m_j, (0.05 \times m_j)^2)$ \tiny /*incorporate uncertainty in $m_j$*/\\
		\normalsize\STATE $c_j = l(m'_j)$ \tiny/*generate measurement value corresponding to $m'_j$ from reference function*/\\
		\normalsize\STATE $c_j' \sim N(c_j, (0.05 \times c_j)^2)$ \tiny/*incorporate uncertainty in $c_j$*/\\
		\ENDFOR
		\normalsize\STATE $g'_i = pfit(m', c')$ \tiny/*fit polynomial quadratic function*/\\
		\normalsize/*begin measurement process*/\\
		\normalsize\FOR{$i = 1-N$}
		\normalsize\STATE $c'_{t(i)} \sim N(c_t, (0.05 \times c_t)^2)$ \tiny/*incorporate uncertainty in $c_{t(i)}$*/\\
		\normalsize\STATE $s'_{t(i)} = g'_i(c'_{t(i)})$ \tiny/*calculate measurand value from $c'_{t(i)}$ and $g_i'$*/\\
		\ENDFOR
		\STATE $u_s = SD(s'_t)$  
	\end{algorithmic}
\end{algorithm}

In the above algorithm, $pfit$ and $SD$ refer to the nonlinear regression method used to characterize the quadratic calibration function and the function used in the program to calculate the standard deviation of the vector of test sample measurements. $m', c', g', c'_t$ and $s'_t$ represent versions of $m, c, g, c_t$ and $s_t$ with uncertainty incorporated into them.

The optimization problem in this situation then involves optimal selection of the calibrators $m = \{m_1, m_2,...,m_n\}$ such that $u_s(m)$ is minimized, subject to constraints on the minimum differences between adjacent calibrator values and integrality restrictions on the calibrator values (manufacturers may have samples with a limited number of calibrator concentration levels they can effectively manufacture). Thus if each calibrator $m_i \in \mathbb{S}_i$, where $\mathbb{S}_i$ is finite and countable, then the feasible region becomes $\prod\limits_{i=1}^{n}\mathbb{S}_i$. The optimization problem thus can be stated as follows.

\begin{equation}
	\label{eq6}
	\begin{aligned}
		&\underset{m \in \prod\limits_{i=1}^{n}\mathbb{S}_i}{\text{min}} \; E[u_s(m)] \simeq \frac{1}{T}\sqrt{\frac{\sum\limits_{i = 1}^N (s_i - \bar{s})^2}{N-1}}\\
		&\text{subject to:}\; \\
		&\frac{m_{i+1} - m_{i}}{m_i} \geq 0.1, \; \forall \; i = 1\; \text{to}\; n-1\\
	\end{aligned}
\end{equation}

In formulation~\ref{eq6} above, $T$ represents the number of replicate measurements of the measurement uncertainty of the assay (that is, $T$ executions of Algorithm~\ref{alg1}). The neighborhood structure for this problem is created as follows. For each calibrator $m_i$, the neighborhood point consists of four points within a distance of 2 units from $m_i$. For example, if the current value of one of the calibrators is 40 units, then its neighborhood is defined as consisting of $N_i(m_i = 40) = \{38,39,41,42\}$ with equal probability of selection of each of the points, as long as it does not violate the minimum distance constraints defined in formulation~\ref{eq6}. Similar neighborhood structures are constructed for each of the six calibrators. In case a new neighborhood point satisfying the constraints is not found for a given calibrator, then it stays at its current value, and the neighborhoods for the other calibrators are explored. The termination criteria for the algorithms is set as achieving a minimum decrease (minimum decrease in Table~\ref{tab3}) in the objective function value within 60 iterations (that is, evaluations of 60 candidate solutions). If the algorithm does not achieve the minimum decrease, for example 5\%, in the value of the objective function within 60 iterations, then this is counted as a failure of the algorithm. The results of applying the original SR method and our proposed relaxation to solve problem~\ref{eq6} is shown in Table~\ref{tab3} below ($\alpha = 0.75$). Note that the stochastic ruler $\theta(a,b)$ was parameterized as $\theta(0,2)$, based on the range of values of $u_s(m)$ observed in preliminary computational experiments.

\begin{table}[htbp]
	\centering
	\caption{Comparison between original and relaxation of the SR method: optimal assay calibrator selection problem.}
	\begin{threeparttable}
	\begin{tabular}{|c|c|c|c|c|c|c|}
		\hline
		& \multicolumn{3}{c|}{Modified SR Method} & \multicolumn{3}{c|}{Original SR Method} \bigstrut\\
		\cline{2-7}    \multicolumn{1}{|C{2cm}|}{Minimum decrease} & \multicolumn{1}{C{1.4cm}|}{Runtime (seconds)} & \multicolumn{1}{C{1.2cm}|}{Average $k$} & \multicolumn{1}{C{1.4cm}|}{Number of failures} & \multicolumn{1}{C{1.4cm}|}{Runtime (seconds)} & \multicolumn{1}{C{1.2cm}|}{Average $k$} & \multicolumn{1}{C{1.4cm}|}{Number of failures} \bigstrut\\
		\hline
		0.050  & 10.672 & 12.607 & 22 & 8.467 & 14.677 & 33 \bigstrut\\
		\hline
		0.075 & 14.425 & 16.091 & 17 & 10.878 & 15.579 & 31 \bigstrut\\
		\hline
		0.100   & 21.293 & 24.269 & 24 & 9.054 & 11.786 & 36 \bigstrut\\
		\hline
		0.125 & 25.024 & 29.417 & 26 & 16.127 & 23.294 & 33 \bigstrut\\
		\hline
	\end{tabular}%
\begin{tablenotes}
	\footnotesize
	\item Minimum decrease: minimum decrease in the objective function value to be achieved prior to termination of the method. 
\end{tablenotes}
\end{threeparttable}
	\label{tab3}%
\end{table}%

Note that the average runtimes and the number of iterations are calculated only for cases where termination criteria were successfully met - that is, only when the minimum decrease was achieved within the computational budget. It is evident that while runtimes are higher for the proposed relaxation of the SR method, the number of failures for the original method is substantially higher for every value of the minimum decrease limit. As expected, the number of failures increase substantially as the minimum decrease limit increases. 

From Table~\ref{tab3}, it appears that when the dimensionality of the solution space increases, as in the case of the optimal calibrator selection problem, the probabilities of failure appear to increase for the SR method, and in such a situation, the proposed relaxation appears to perform better with significantly fewer failures. In our final numerical example, we consider another multi-dimensional discrete simulation optimization problem within the context of facility location, and focus on comparing the number of failures and the optimal objective function values for the proposed relaxation and the original version of the SR method.

\subsection{Numerical Example 4: Optimal Facility Location}
\begin{figure}
	\centering
	\subfloat[Neighbors of $x$ not located on edges or corners of the grid.]{%
		\resizebox*{5cm}{!}{\includegraphics{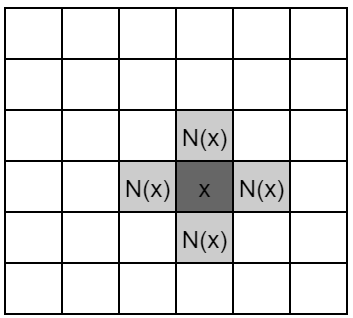}}}\hspace{30pt}
	\subfloat[Neighbors of $x$ located on edges or corners of the grid.]{%
		\resizebox*{5cm}{!}{\includegraphics{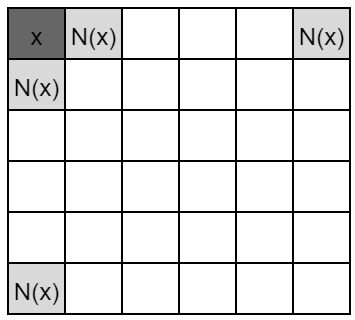}}}
	\caption{Neighborhood structure for the facility location simulation optimization problem. $x$ represents the current location of a facility, and the $N(x)$ represent the neighbors of $x$.} \label{fig:facloc}
\end{figure}

We now construct another numerical example wherein the global optima are not known - a stochastic facility location problem with a discrete solution space. Our formulation of this facility location problem can be considered to be a discretization of a continuous facility location problem. This is reflected in our formulation of the solution space: it is a 6 $\times$ 6 rectangular grid, where each grid square is a potential solution. This is depicted in Figure~\ref{fig:facloc}, where we discuss the neighborhood structure that we impose on the discrete solution space. 

We consider the location of three facilities within this rectangular grid, implying that there are three decision variables. We introduce stochasticity in this problem by modelling the daily demand from each grid square (to be serviced by one of the three facilities), independent of the demand from other grid squares, as a Gaussian random variable, with parameters $N(180, 30^2)$. This parameterization of the daily demand implies that negative values of demand are highly unlikely to occur, and can hence reasonably be ignored (and demand can be resampled) when they occur. We assume that the demand from a grid square is serviced by the facility it is closest to, with distance calculated rectilinearly. The objective of the optimization problem thus becomes locating three facilities in the grid such that the average distance travelled to the facilities from across the entire grid is minimized. We calculate the average distance travelled to these facilities over some pre-determined time period. The simulation conceptualization of this facility location problem is described in Algorithm~\ref{alg:fac} below. 

\begin{algorithm}
	\caption{Simulation conceptualization of the facility location problem.}
	\label{alg:fac}
	\begin{algorithmic}
		\STATE Initialize with rectangular grid, grid square index $j \in \{1,2,..,J\}$, simulation duration $T_0$ days, daily demand distribution $N(\mu, \sigma^2)$. 
		\FOR{$t = 1$ to $T_0$}
		\FOR{$j = 1$ to $J$}
		\STATE Sample demand of $j^{th}$ grid square $d_j \sim N(\mu, \sigma^2)$
		\STATE Calculate distance travelled to nearest facility to satisfy demand $d_j$ as $dist_j$
		\STATE Calculate total distance travelled from square $j$ on $t^{th}$ day as $tot_{jt} = d_j \times dist_j$
		\ENDFOR
		\STATE Calculate average distance travelled across grid on $t^{th}$ day as $dist_t = \frac{\sum\limits_{j = 1}^{J} tot_{jt}}{J}$
		\ENDFOR
		\STATE Calculate average daily distance travelled across grid as $dist = \frac{\sum\limits_{t = 1}^{T_0} dist_t}{T_0}$
	\end{algorithmic}
\end{algorithm}

Thus a single replication of the simulation described in Algorithm~\ref{alg:fac} yields the value $dist$, and hence $dist$ represents $h(x)$, where $x$ is a three-tuple of locations representing the decision variables. We specify the neighborhood structure for solutions $x$ by doing so for each of its components - i.e., by specifying a neighborhood structure for each candidate location for a facility. For this purpose, we define the coordinates of a grid as $i, j$, where both $i, j \in \{1,2,3,4,5,6\}$ and $(i,j) = (1,1)$ corresponds to the bottom left-most grid square. For the candidate location $(i, j)$ where $1 < i, j < 6$, we define its neighbours as follows $N(i,j) = \{(i-1, j), (i+1, j), (i, j-1), (i, j+1)\}$. For squares located at the corners or the edges of the grid - i.e., for $(i,j)$ such that $i \in \{1,6\}$ and/or $j \in \{1,6\}$ - the neighbors are defined as follows. For $(i, j) = (1, j)$ where $1 < j < 6$, $N(1,j) = \{(6, j), (2, j), (1, j - 1), (1, j + 1)\}$. For $(i, j) = (i, 1)$ where $1 < i < 6$, $N(i,1) = \{(i, 6), (i, 2), (i - 1, 1), (i + 1, 1)\}$. The neighborhood of a given solution $x$ is then specified as the Cartesian product of the neighborhoods of its components, where the neighborhoods of each component are defined as above.

%

The next step involved parameterizing the stochastic ruler $\theta(a,b)$ - i.e., estimating appropriate values of $a$ and $b$. Given that the global minima for this problem are not known, we estimated these parameters via extreme case calculations as $a = 250$ and $b = 800$. We then instituted the following the termination criteria for the algorithm: a prespecified percentage reduction in the objective function value is achieved (ranging from 5\% to 50\% in our numerical experiments), as estimated by the mean of the samples generated (in Step 2 of the stochastic ruler method) at a given candidate solution, or when the total number of iterations (i.e., the value of $k$ in either the original or modified stochastic ruler method) reaches 150. We now provide the comparison between the original version and our proposed modification of the stochastic ruler method in Table~\ref{tab:fac}. 

\begin{table}[htbp]
	\centering
	\caption{Comparison between the original version and relaxation of the SR algorithm: optimal facility location problem.}
	\begin{threeparttable}
		\begin{tabular}{|c|c|c|c|c|}
			\hline
			\multicolumn{1}{|c|}{\multirow{2}[4]{*}{Minimum decrease}} & \multicolumn{2}{m{10.045em}|}{Number of failures} & \multicolumn{2}{c|}{$\bar{h}(x_*)$} \bigstrut\\
			\cline{2-5}          & \multicolumn{1}{p{6em}|}{Relaxation} & \multicolumn{1}{p{4.045em}|}{Original} & \multicolumn{1}{p{5.865em}|}{Relaxation} & \multicolumn{1}{p{4.045em}|}{Original} \bigstrut\\
			\hline
			0.05  & 1     & 33    & 501.2 & 521.6 \bigstrut\\
			\hline
			0.10   & 1     & 33    & 501.2 & 521.6 \bigstrut\\
			\hline
			0.15  & 1     & 34    & 481.4 & 504.3 \bigstrut\\
			\hline
			0.25  & 0     & 46    & 462.6 & 497.4 \bigstrut\\
			\hline
			0.50   & 0     & 76    & 383.1 & 387.8 \bigstrut\\
			\hline
		\end{tabular}%
		\begin{tablenotes}
			\footnotesize
			\item Minimum decrease: minimum decrease in the objective function value to be achieved prior to termination of the method. Original: original version of the SR method. Relaxation: our proposed relaxation of the SR method. $\bar{h}(x_*)$ = average objective function value at the optimal solution $x_*$. 
		\end{tablenotes}
	\end{threeparttable}
	\label{tab:fac}%
\end{table}%

As discussed earlier, for this numerical example, we focus on the number of failures and the objective function value for comparing the performance between the original version and the proposed relaxation of the SR method. It is evident from Table~\ref{tab:fac} that the original version of the SR method exhausts its computational budget a significantly larger number of times than our proposed relaxation without achieving the minimum decrease required in the objective function value. 

\section{Discussion}
\label{disc}
In this paper, we present preliminary results regarding a relaxation to the stochastic ruler method originally developed by \cite{yan1992}. We have provided an overview of how the proposed relaxation fits within the framework for convergence to the global optimum developed by Yan and Mukai. We have also provided detailed numerical examples to demonstrate the performance of our proposed relaxation to the original SR method. In the first two numerical examples, where the global minimum is known \textit{a priori}, we show that our proposed relaxation outperforms the original method. In both these examples, the dimensionality of the search space is one, and the neighborhood structures are simple. In the third numerical example, which is constructed on the basis of the real-world problem of optimal calibrator selection for clinical measurement processes \cite{wsc18vr}, we see that the original SR method appears to have higher runtimes, but at the expense of significantly higher failures (i.e., when algorithm exceeds the computational budget before achieving the minimum decrease in the objective function value). This implies that if the runtimes corresponding to the failures are considered, the average computational runtimes for the relaxation are lower than that for the original SR method. In the fourth numerical example, which again involves a multi-dimensional solution space and neighborhood structure, we see that the number of failures for the original SR method are substantially higher, and that the objective function values achieved at the solutions identified by the relaxation are in all cases lower than that achieved by the original version. These numerical examples illustrate the key advantage of our proposed relaxation, in that it decreases the computational expense required for selecting promising candidate solutions as the next estimate of the optimal solution.

Several avenues of future research suggest themselves, a few of which we are currently pursuing. First, we are exploring whether similar relaxations and their theoretical underpinnings can be developed for the modifications of the SR method presented in \cite{alrefaei2001mod,alrefaei2005}. Indeed, our preliminary numerical experiments indicate that such relaxations are possible. A second possibility involves an adaptive version of the SR method that can be developed - for example, the value of $\alpha$ can change as a function of $k$ as the algorithm converges towards the global optimum. Our numerical examples also illustrate the fact that further investigation regarding the appropriate choice of $\alpha$ is also required.

\bibliographystyle{apacite}
\bibliography{srbib}

\end{document}